\documentclass{amsart}
\usepackage{amsmath,amsfonts,amscd,color}
\usepackage{latexsym}
\usepackage[latin1]{inputenc}
\definecolor{orange}{rgb}{1,0.5,0}
\title[Index theory of boundary value problems]{C$^*$-Algebra approach to the index theory of boundary value problems}
\author[S. Melo, T. Schick, E. Schrohe]{Severino T. Melo, Thomas Schick \and Elmar Schrohe}

\newtheorem{thm}{Theorem}
\newtheorem{pro}[thm]{Proposition}
\newtheorem{lem}[thm]{Lemma}
\newtheorem{cor}[thm]{Corollary}
\newtheorem{df}[thm]{Definition}
\newtheorem{rem}[thm]{Remark}

\newcommand{\asim}{\!\!\sim}
\newcommand{\A}{{\mathcal A}}
\newcommand{\ac}{{\mathfrak A}}

\newcommand{\is}{{\mathcal I}}

\newcommand{\ks}{{\mathcal K}}
\newcommand{\ts}{{\mathcal T}}
\newcommand{\ws}{{\mathcal W}}
\newcommand{\kc}{{\mathfrak K}}
\newcommand{\hc}{{\mathfrak H}}
\newcommand{\R}{{\mathbb R}}

\newcommand{\C}{{\mathbb C}}

\newcommand{\cA}{{\mathcal A}}

\newcommand{\lh}{{\mathcal L}(\mathcal H)}
\newcommand{\naturals}{{\mathbb N}}
\newcommand{\Z}{{\mathbb Z}}

\newcommand{\diff}{\mbox{Diff}}
\newcommand{\ind}{\mbox{\rm ind}\,}
\newcommand{\op}{operator}
\newcommand{\ops}{operators}

\newcommand{\End}{\mbox{{\sc End}}}
\newcommand{\cqd}{\hfill$\Box$}
\newcommand{\h}{{\mathcal H}}

\newcommand{\intX}{X^\circ}
\newcommand{\rp}{{\mathbb R}_{\ge0}}
\newcommand{\iso}{\cong}
\newcommand{\complexs}{\mathbb{C}}
\DeclareMathOperator{\id}{id}

\begin{document}

\begin{abstract}
Boutet de Monvel's calculus provides a pseudodifferential framework which encompasses the classical differential boundary value problems. In an extension of the concept of Lopatinski and Shapiro, it associates 
to each operator two  symbols: a pseudodifferential principal symbol, which is
a bundle homomorphism, and  an operator-valued boundary symbol. Ellipticity
requires the invertibility of both. If the underlying manifold is compact,
elliptic elements define Fredholm operators. Boutet de Monvel \cite{B} showed
how then the index can be computed in topological terms. The crucial
observation is that elliptic operators can be mapped to compactly supported
$K$-theory classes on the cotangent bundle over the interior of the
manifold. The Atiyah-Singer topological index map, applied to this class, then
furnishes the index of the operator. Based on this result, Fedosov,
Rempel-Schulze and Grubb have given index formulas in terms of the symbols. In
this paper we survey how C$^*$-algebra K-theory, as
initiated in \cite{MNS}, can be used to give a proof of Boutet de Monvel's 
index theorem for boundary value problems, a task carried out in \cite{MSS}, 
and how the same techniques 
yield an index theorem for families of Boutet de Monvel operators, detailed
in \cite{MSSf}.  The key ingredient of our 
approach is a precise description of the K-theory of the kernel and of the
image of the boundary symbol.
\end{abstract}
\maketitle
\begin{center}{\footnotesize 2010 Mathematics Subject Classification: 19K56, 46L80, 58J32}
\end{center}

\section{Boutet de Monvel's calculus}
Let $X$ be a compact $n$-dimensional manifold with boundary $\partial X$, 
embedded in a closed manifold $\widetilde X$ of the same dimension.
By $\intX$ we denote the interior of $X$.
We assume that $X$ is connected and $\partial X$ is nonempty.
Given a pseudodifferential operator $P$ on $\widetilde X$, we define
the truncated pseu\-do\-dif\-ferential operator 
$P_+\colon C^\infty(X)\to C^\infty(\intX)$ 
as the composition
$r^+Pe^+$, where $e^+$ is extension by zero from $X$ to $\widetilde X$ and 
$r^+$ is the restriction to $\intX$. In general, the functions in the 
range of $P_+$ will not be smooth up to the boundary. 
One therefore assumes that $P$ satisfies the {\em transmission
condition}, a condition on the symbol of $P$ which we recall
in \eqref{eq:transmit_prop} and
which ensures that both $P_+$ and $(P^*)_+$, the truncated 
operator formed from the formal adjoint of $P$, map smooth 
functions on $X$ to smooth functions on $X$. 

An operator in Boutet de Monvel's calculus is a matrix
\begin{eqnarray}
\label{1}
A=
\begin{pmatrix}P_++G&K\\T&S
\end{pmatrix}: 
\begin{array}{ccc}
C^\infty(X,E_1)& & C^\infty(X,E_2)\\
\oplus&\to&\oplus\\
C^\infty(\partial X,F_1)&&C^\infty(\partial X,F_2)
\end{array} 
\end{eqnarray} 
acting on sections of vector bundles $E_1$, $E_2$ over $X$  and $F_1$, $F_2$ over $\partial X$.
Here, $P$ is a pseudodifferential operator satisfying the transmission condition; 
$G$ is a singular Green operator, 
$T$ is a trace operator, 
$K$ is a potential (or Poisson) operator, 
and $S$ is a pseudodifferential operator on $\partial X$. 
All these operators are assumed to be classical; i.e.\ their symbols have polyhomogeneous expansions 
in the respective classes.
The calculus contains the classical boundary value problems, 
where $P$ is a differential operator, $G=0$, and $T$ a differential trace operator.
Here $F_1=0$; the operators $K$ and $S$ do not appear.   
It also contains their inverses, provided they exist. 
In this case, $F_2=0$, the operators $T$ and $S$ do not show up,  
and the inverse to $\binom {P_+}T$ is of the form $(Q_++G\ \ K)$, 
where  $K$ solves
the semi-homogeneous problem $Pu=0, Tu=g$ for given $g$, and $Q_++G$ solves the 
semi-homogeneous problem $Pu=f, Tu=0$ for given $f$. 
Here $Q$ is a parametrix to $P$, and $G$ is the correction needed to fulfill the boundary condition.
For details, we refer to the monographs
by Rempel and Schulze \cite{RS} or Grubb \cite{G} as well as to the short introduction \cite{S3}.

The operators $G$, $K$, and $T$\
are regularizing in the interior of $X$. 
In a collar neighborhood of the boundary, they can be viewed as 
operator-valued pseudodifferential operators along the boundary. 
In particular, they have an order assigned to them.
The singular Green and the trace operators 
also have a {\em class}\ (or {\em type}) $d\in\naturals_0$, 
related to the order of the derivatives appearing in the 
boundary condition. 
The composition of two operators  of the form \eqref{1}
is defined whenever the vector bundles serving as the range of the first operator 
form  the domain of the second.
The composition $AA'$ of an operator $A'$ of order $m'$ and class 
$d'$ with an operator $A$ of order $m$ and
class $d$ results in an operator of order $m+m'$ and  class $\le \max(m'+d, d')$. 
In particular, the composition of two operators of order and 
class zero is again of order and class zero.

For $E_1=E_2=E$ and $F_1=F_2=F$, the operators of order and class 
zero thus form an algebra $\cA^\circ$. 
Moreover, they extend to bounded operators on the 
Hilbert space $\h=L^2(X,E)\oplus L^{2}(\partial X,F)$.
In fact, $\cA^\circ$ is a $*$-subalgebra of the algebra $\lh$ of all bounded operators on $\h$,
closed under holomorphic functional calculus, cf.\ \cite{S2}.%
\footnote{We use here the definition of order and class in \cite{RS} and \cite{S3} 
 which differs slightly from that in \cite{G}.
 It allows us to use the $L^2$-space over the boundary instead of the Sobolev space $H^{-1/2}(\partial X,F)$
 and gives us better homogeneity properties of the boundary symbols.
 As the kernel and the cokernel of an elliptic 
 operator in $\cA^\circ$ consist of smooth functions, the choice is irrelevant 
 for index theory.}

Standard reductions --recalled in \cite[Section 1.1]{MSS}--
allow to reduce an arbitrary  index problem in the calculus (defined by an elliptic Boutet de Monvel operator of arbitrary order and class and acting between different bundles)
to the case where the order and class are zero and  $E_1=E_2=E$ and $F_1=F_2=F$. In other words, it suffices to study the index problem for ellipitic operators in $\cA^\circ$ and we are then free to apply operator-algebraic methods. There is also no loss of generality in the assumption that the manifold $X$ is connected.

We consider an operator $A$ as in  \eqref{1}, with 
$E_1=E_2=E$, $F_1=F_2=F$. 
The pseudodifferential principal symbol $\sigma(A)$ of $A$ is defined to be the 
principal symbol of the pseudodifferential part $P$ (a smooth bundle morphism), restricted to $S^*X$. 
This makes sense as $G$ is regularizing in the interior. The choice of a hermitian structure on $E$ (already needed to define the inner-product of $\h$) turns the map 
\[
\cA^\circ\ni A\mapsto \sigma(A)\in\mbox{{\sc HOM}}(\pi^*E) 
\] 
into a homomorphism of $*$-algebras. We have denoted by $\pi\colon S^*X\to X$ the canonical projection of the co-sphere bundle of $X$. 

The boundary principal symbol of an $A\in\cA_0$ is a smooth endomorphism of 
\begin{equation}\label{end}
\left( L^2(\rp)\otimes \pi_\partial^*E|_{\partial X}\right)\oplus\pi_\partial^* F, 
\end{equation}
with $\pi_\partial\colon S^*\partial X\to \partial X$ denoting the canonical projection of the co-sphere bundle of $\partial X$.
It is best described for a trivial one-dimensional bundle and 
in local coordinates $(x',x_n,\xi',\xi_n)$ for 
$T^*X$ in a neighborhood of the boundary.
Here, $G$ acts like a pseudodifferential operator along
the boundary, with an operator-valued symbol taking values 
in regularizing operators in the normal direction.
One way to write this operator-valued symbol is via a so-called 
symbol kernel
$\tilde g = \tilde g(x',\xi',x_n,y_n)$. 
For fixed $(x',\xi')$, this is a rapidly decreasing function 
in $x_n$ and $y_n$ which acts as an integral operator on $L^2(\rp)$. 
It satisfies 
special estimates, combining the usual pseudodifferential
estimates in $x'$ and $\xi'$ with those for rapidly decreasing functions in
$x_n$ and $y_n$. 
The singular Green symbol $g$ of $G$ is defined from the symbol kernel via 
Fourier and inverse Fourier transform:
$$
g(x',\xi',\xi_n,\eta_n) = F_{x_n\to \xi_n}\overline F_{y_n\to \eta_n}\tilde
g(x',\xi',x_n,y_n).
$$
It has an expansion into homogeneous terms; the leading one we call $g_0$. Inverting
the operation above, we associate with $g_0$\ a symbol kernel 
$\tilde g_0(x',\xi',x_n,y_n)$\ which is rapidly decreasing in $x_n$ and $y_n$
for fixed $(x',\xi')$. 
We denote by $g_0(x',\xi',D_n)$ the (compact) operator
induced on $L^2(\rp)$ 
by this kernel. 
Similarly,  $K$ and $T$ have symbol-kernels $\tilde k(x',\xi',x_n)$ and 
$\tilde t(x',\xi',y_n)$; these are  rapidly decreasing 
functions for fixed $(x',\xi')$. 
The symbols $k$ and $t$ are defined as their Fourier and inverse 
Fourier transforms. They have asymptotic expansions with leading terms
$k_0$ and $t_0$. Via the symbol-kernels $\tilde k_0$ and $\tilde t_0$
one defines $k_0(x',\xi', D_n)\colon \complexs\to L^2(\rp)$  
as multiplication by $\tilde k_0(x',\xi',\cdot)$, while 
$t_0(x',\xi',D_n) :L^2(\rp)\to \complexs$ 
is the
operator $\varphi \mapsto \int \tilde t_0(x',\xi',y_n)\varphi(y_n)\,dy_n$.   

We denote by $p_0$ and $s_0$ the principal symbols of $P$ and $S$, respectively.
The boundary symbol $\gamma(A)$ of $A$ at $(x',\xi')$ is then defined by 
\begin{equation}\label{bdysymbol}
\gamma(A)(x',\xi')=
\left(\begin{array}{cc}
p_0(x^\prime,0,\xi^\prime,D_n)_{_{+}}+g_0(x',\xi',D_n)&k_0(x',\xi',D_n)\\
t_0(x',\xi',D_n)&s_0(x',\xi')
\end{array}\right).
\end{equation}

This gives an invariantly defined operator-valued function on $T^*\partial X$ only up to a choice of a normal coordinate; i.e., we need to restrict ourselves to an atlas whose changes of coordinates, near the boundary, preserve not only the boundary points 
$\{x_n =0\}$ but the variable $x_n$ as well \cite[Theorem 2.4.11]{G}. The boundary symbol can be viewed as a function on $S^*\partial X$ due to its {\em twisted homogeneity}, 
\begin{equation}
\begin{pmatrix}\kappa_{\lambda^{-1}}&0\\0&\id\end{pmatrix}
\gamma(A)(x',\lambda\xi')
\begin{pmatrix}\kappa_{\lambda}&0\\0&\id\end{pmatrix}=\gamma(A)(x',\xi'),\quad\lambda>0,
\end{equation}
with the $L^2(\rp)$-unitary $\kappa_\lambda$ given by $\kappa_\lambda f(t)=\sqrt\lambda f(\lambda t)$.

A connection between Toeplitz operators and pseudodifferential operators satisfying the transmission condition turns out to be an essential point for both the computation of the 
$K$-theory of the range of the principal boundary symbol and for the proof of the estimate (\ref{estimate}) needed to describe its kernel.
Let $p\sim\sum p_j$ be the asymptotic expansion of the local symbol 
$p$ of $P$ into terms $p_j(x,\xi)$, 
which are positively homogeneous of degree $j$ in $\xi$ for $|\xi|\ge 1$. 
The transmission condition requires that, for $x_n=0$ and $\xi=(0,\pm1)$, 
\begin{equation}\label{eq:transmit_prop}
D_{x}^\beta D^\alpha_\xi p_j(x',0, 0,1) = (-1)^{j-|\alpha|}
D_{x}^\beta D^\alpha_\xi p_j(x',0, 0,-1).
\end{equation}
Hence the limits of 
$p_0(x',0,\xi',\xi_n)$ as $\xi_n\to\pm\infty$ coincide for fixed 
$(x',\xi')$, and the function 
$$p_{(x',\xi')}(z)=p_0\left (x',0,\xi',\frac{iz-i}{z+1}\right),\ z\in S^1,z\neq -1,$$
extends continuously to $S^1$.

We next observe that the image of the Hardy space $H^2(S^1)$ under the unitary
map $U\colon L^2(S^1)\to L^2(\R)$,
\[
Ug(t)=\frac{\sqrt{2}}{1+it}g\left(\frac{1-it}{1+it}\right),
\]
is equal to $F(L^2(\R_{\geq 0}))$, where $F$ denotes the Fourier transform. The truncated Fourier multiplier 
\begin{eqnarray}\label{trunc}
\begin{array}{rcl}
p_0(x',0,\xi',D_n)_{_{+}}\colon L^2(\R_{\geq 0})&\longrightarrow&L^2(\R_{\geq 0})\\
u&\longmapsto&F^{-1}(p_0(x',0,\xi',\cdot)Fu)|_{\R_{\geq 0}}
\end{array}
\end{eqnarray}
is therefore equal to $F^{-1}UT_{p_{x',\xi'}}U^{-1}F$, where $T_{p_{x',\xi'}}$ denotes the Toeplitz operator of symbol $p_{x',\xi'}$. It then follows from classical results about Toeplitz operators \cite{Dou}\ that
\begin{equation}
\label{toeplitz}
\|p_0(x',0,\xi',D_n)_{_{+}}\|=\sup_{\xi_n}|p_0(x',0,\xi',\xi_n)|=\inf_{K}\|p_0(x',0,\xi',D_n)_{_{+}}+K\|,
\end{equation}
with the last infimum being taken over all compact operators $K$ on $L^2(\rp)$. In particular, $p_0(x',0,\xi',D_n)_{_{+}}$ is compact if and only if $p_0(x',0,\xi',\xi_n)=0$ for all $\xi_n\in\R$. 

Gohberg \cite{Goh}\ and Seeley \cite{See}\ established the equality between the norm, modulo compacts, of a singular integral operator on a compact manifold and 
the supremum norm of its symbol. Proofs of that estimate in the language of pseudodifferential appeared in \cite{Hp,KN}.  
The following generalization holds for Boutet de Monvel operators (a proof for this result can be found in Rempel and Schulze's book \cite[2.3.4.4]{RS}; they credit Grubb and Geymonat \cite{GG}\ for earlier work):
\begin{equation}
\label{rs}
\inf_{C\in\ks}\|A+C\|=\max\{\|\sigma(A)\|,\|\gamma(A)\|\}, \ \text{for all}\ A\in\cA^\circ,
\end{equation}
with $\ks$\ denoting the ideal of the compact operators on $\h$, $\|\sigma(A)\|$ the supremum norm of $\sigma(A)$\ on $S^*X$, and $\|\gamma(A)\|$ the supremum
over all $(x',\xi')$\ in $S^*\partial X$\ of $\|\gamma(A)(x',\xi')\|$.  

\begin{df}\label{defas}
We denote by $\cA$ the norm closure in $\lh$ of the algebra $\cA^\circ$ of all classical Boutet de Monvel operators of order and class zero.
\end{df}

Equation (\ref{rs}) implies, in particular, that $\sigma$ and $\gamma$ extend to C$^*$-algebra homomorphisms defined on $\cA$ and taking values in continuous endomorphisms of the bundles $\pi^*E$ and 
$\left( L^2(\rp)\otimes \pi_\partial^*E|_{\partial X}\right)\oplus\pi_\partial^* F$, respectively.
It also implies that $\ker\gamma\cap\ker\sigma=\ks$ and that the quotient $\cA/\ks$ is isomorphic to the image of the pair $(\sigma,\gamma)$ and, in particular, that $A\in\cA$ is Fredholm if and only if both $\sigma(A)$ and $\gamma(A)$ are invertible. This description of $\cA/\ks$, however, is not explicit enough for $K$-theory computations. 

In the rest of the paper, we will assume that $E=X\times\C$ and $F=\partial
X\times\C$. For the general case, all the results can be reformulated
in a straightforward way and their proofs can be 
adapted. A key ingredient in this adaption are canonical Morita equivalences
between the algebras acting on 
functions, the algebras acting on vectors of functions and the algebras acting
on sections of general bundles, giving rise to canonical K-theory
isomorphisms. They are based on the well known \cite[Subsection 1.5]{MNS}
Morita equivalence between 
sections of the endomorphism bundle $\End(E)$ of a bundle $E$ and the algebra
of functions itself, given by the bimodule of sections of $E$.

\section{The boundary-symbol exact sequence}

Our description of the kernel of the boundary symbol (or rather of its
quotient by the compacts) depends on an estimate for the norm, modulo an ideal
of operators in $\cA^\circ$ whose closure is larger than the compacts. 

\begin{thm} \label{kernelgamma}
The principal symbol $\sigma$ induces a C$^*$-algebra isomorphism 
\[
\ker\gamma/\ks
\ni[A]\longmapsto \sigma(A)\in C_0(S^*X^\circ),
\]
where $C_0(S^*X^\circ)$ denotes the algebra of continuous functions on $S^*X$ which vanish at the boundary.
\end{thm}
{\em Sketch of proof}: If the upper left corner of the matrix in the right 
hand side of (\ref{bdysymbol}), for an $A\in\cA^\circ$, vanishes, then $p_0(x',0,\xi',D_n)_{_{+}}$ 
is compact, since $g_0(x',\xi',D_n)$ is compact. 
It then follows from (\ref{toeplitz}) that $p_0(x',0,\xi',\xi_n)=0$ for all $\xi_n\in\R$ and hence $g_0(x',\xi',D_n)=0$. 
This shows that the kernel of $\gamma$ restricted to $\cA^\circ$ is equal to the set $\is^\circ$ of all $A$ as in (\ref{1}) such that $\sigma(A)$ vanishes at the boundary and, 
moreover, $G$, $K$, $T$ and $S$ are of lower order. To prove that the kernel of $\gamma$ (defined on the whole algebra $\cA$) is equal to the closure of 
$\is^\circ$, which we denote by $\is$, one needs to use that there exists $C>0$ such that
\begin{equation}\label{estimate}
\inf\{\|A+A'\|, A'\in\is^\circ\}\leq C \|\gamma(A)\|
\end{equation} 
for all $A\in\cA^\circ$. The proof of this estimate \cite[Lemma 2]{MNS} combines the above mentioned classical Gohberg-Seeley estimate with 
(\ref{toeplitz}).

The closed ideal $\is$ can also be described as the C$^*$-subalgebra of $\lh$ generated by all the operators of the form
\begin{eqnarray}\label{defis}
A=
\begin{pmatrix}\varphi P\varphi &K_1\\K_2&K_3
\end{pmatrix}: 
\begin{array}{ccc}
L^2(X)& & L^2(X)\\
\oplus&\to&\oplus\\
L^2(\partial X)&&L^2(\partial X)
\end{array} ,
\end{eqnarray} 
where $P$ is a zero-order classical pseudodifferential operator, $\varphi$ is (multiplication by) a smooth function with support contained in $X^\circ$, and $K_1$, $K_2$ and $K_3$ are compact operators.
It then follows from the Gohberg-Seeley estimate that the principal symbol induces the desired isomorphism.
\cqd.

Given $f\in C(X)$, the operator $m(f)$ defined by
\[
L^2(X)\oplus L^2(\partial X)\ni (\phi,\psi)\mapsto m(f)(\phi,\psi)=(f\phi,0)\in L^2(X)\oplus L^2(\partial X)
\]
belongs to $\cA$. Abusing notation a little, let us denote also by $m\colon C(X)\to\cA/\ks$ the C$^*$-algebra homomorphism that maps $f$ to the class of $m(f)$ in the quotient $\cA/\ks$ and also by $\gamma$ the map induced by the boundary symbol on the quotient $\cA/\ks$ with kernel $\ker\gamma/\ks$. Taking into account that the isomorphism of Theorem~\ref{kernelgamma} is induced by the principal symbol and that the principal symbol of the multiplication by a function is equal to the function itself, we then get the following commutative diagram of C$^*$-algebra exact sequences
\begin{equation}
\label{commdiag}
\def\mapup#1{\Big\uparrow\rlap{$\vcenter{\hbox{$\scriptstyle#1$}}$}}
\begin{array}{ccccccc}
0\longrightarrow&C_0(S^*X^\circ)&\longrightarrow&\cA/\ks&{\mathop{\longrightarrow}\limits^{\gamma}}&\mbox{im}\,\gamma&\longrightarrow 0
\\                &\mapup{m^\circ}&&\mapup{m}&&\mapup{b}&
\\0\longrightarrow&C_0(X^\circ)&\longrightarrow&C(X)&{\mathop{\longrightarrow}\limits^{r}}&C(\partial X)&\longrightarrow 0
\end{array},    
\end{equation}
where $m^\circ$ denotes composition with the bundle projection, $r$ denotes the restriction map and $b$ denotes the homomorphism that maps a function $g\in C(\partial X)$ to the boundary principal symbol of $m(f)$ for some $f\in C(X)$ such that $g=r(f)$.

Let $\ts$  denote the Toeplitz algebra on $S^1$. 
It is well-known that $\ts$ contains the compact 
operators and that, as a C$^*$-algebra, $\ts$ is generated by 
the operators $T_\varphi$ for $\varphi\in C(S^1)$. By $\ts_0$ we denote the ideal in $\ts$
generated by the operators $T_\varphi$ with $\varphi$ vanishing at $-1$.

By $\ws^{1,1}$ we denote the image of the Toeplitz algebra in  
$\mathcal{L}(L^2(\R_{\geq 0}))$ under the isomorphism sketched after \eqref{trunc}, i.e.,
$\ws^{1,1}$ is the C$^*$-algebra generated by the truncated Fourier
multipliers  $\varphi(D)_+$, where $\varphi\in C(\mathbb R)$ has equal limits at $\pm\infty$.
We write $\ws^{1,1}_0$ for the corresponding image of $\ts_0$. 

Next, we let $\ws $ denote the C$^*$-subalgebra of 
$\mathcal{L}(L^2(\R_{\geq 0})\oplus \C)$
consisting of all elements whose upper left corner belongs to $\ws^{1,1}$, and by $\ws_0$
the ideal where the upper left corner is in $\ws^{1,1}_0$. $\ws$ is the algebra 
of {\em Wiener-Hopf} operators on $\R_{\ge0}$.
The following observation will play an important role: 
\begin{lem}\label{kw}
We have $K_0(\ws_0)=0=K_1(\ws_0)$. 
\end{lem}
\proof\ Denote, for the moment,  by $\ks$ and $\ks_\oplus$ the compact operators 
on $L^2(S^1)$ and on $L^2(\R_{\geq 0})\oplus \C$, respectively. Then 
$C(S^1)\cong\ts/\ks\cong \ws/\ks_\oplus$, and we have a short exact sequence
$$0\to \ks_\oplus\to \ws\to C(S^1)\to 0,$$
where the map $\ws\to C(S^1)$ is induced by 
$\varphi(D)_+\mapsto \varphi(\frac{iz-i}{z+1}).$  
The associated 6-term exact sequence is 
\begin{equation}
\label{T6}
\begin{array}{ccccc}
\Z&\longrightarrow&K_0(\ws)&\longrightarrow&\Z\\
\big\uparrow&&&&\big\downarrow\\\Z&\longleftarrow&K_1(\ws)&\longleftarrow& 0.
\end{array}
\end{equation}
As there exists a Toeplitz \op\ of Fredholm index one, there also exists an 
\op\ in $\ws$\ of index one;
hence, the index mapping in (\ref{T6}) is surjective. This gives 
$
K_0(\ws)=[I]\cdot\Z\ \text{and}\ K_1(\ws)=0.
$
The six-term exact sequence associated to  
\[
0\longrightarrow\ws_0\longrightarrow\ws\smash{\mathop{\longrightarrow}}
\ \C\longrightarrow 0 
\] 
then shows that $K_{0}(\ws_0)=K_{1}(\ws_0)=0$.\cqd

The rest of this section is devoted to a sketch of the proof of:
\begin{thm} \label{thmimg} The injective C$^*$-algebra homomorphism $b:C(\partial X)\to \mbox{{\em Im}}\,\gamma$ induces a $K$-theory isomorphism.\end{thm}

It follows from our remarks preceding (\ref{toeplitz}) that the image of $\gamma$ is contained in $C(S^*\partial X,\mathcal{W})$.
%
Using standard arguments of the Boutet de Monvel calculus, one shows 
that $C(S^*\partial X,\mathcal{W}_0)$ is contained in the image of $\gamma$, see
\cite[Section 3]{MNS} for details.

Since the intersection of $\mbox{im}\,b$ and $C(S^*\partial X,\mathcal{W}_0)$ is trivial, we have 
\begin{equation}\label{img}
C(S^*\partial X,\mathcal{W}_0)\oplus\mbox{im}\,b\subseteq\mbox{im}\,\gamma,
\end{equation}
where $\oplus$ denotes the direct sum of Banach spaces, not of C$^*$-algebras.
To prove that the reverse inclusion also holds, we need to consider the C$^*$-algebra homomorphism $\lambda$ of $C(S^*\partial X)\otimes\mathcal{W}$ into itself defined by
\[
f\otimes\left(\begin{array}{cc}p(D)_+&*\\ *&*\end{array}\right)\mapsto 
p(\infty)f\otimes \left(\begin{array}{cc}\mbox{Id}&0\\ 0&0\end{array}\right),
\] 
where $\mbox{Id}$ denotes the identity operator on $L^2(\R_{\geq 0})$. If $F\in\mbox{im}\,\gamma$, then $F-\lambda(F)$ belongs to 
$C(S^*\partial X,\mathcal{W}_0)$, and hence to $\mbox{im}\,\gamma$, and then also $\lambda(F)$ belongs to $\mbox{im}\,\gamma$.

Next denote by $\gamma_{11}$ the upper-left corner of $\gamma$ and
suppose that $f\in C(S^*\partial X)$ is such that $(x',\xi')\mapsto f(x',\xi')\mbox{Id}$ belongs 
to the image of $\gamma_{11}$. 
Given $\varepsilon>0$  there exist a pseudodifferential operator  
$P$ with principal symbol $p_0$
and a singular Green operator $G$ with principal symbol $g_0$ 
such that for all $(x',\xi')$ in $S^*\partial X$ 
\begin{eqnarray}\label{est1}
\lefteqn{\|p_0(x',0,\xi',D_n)-g_0(x',\xi',D_n) -f(x',\xi')\mbox{Id}\|}\nonumber\\
&=&
\|\gamma_{11}(P_++G)-f\otimes \mbox{Id}\|<\varepsilon
\end{eqnarray}
As $g_0(x',\xi',D_n)$ is compact, we conclude from \eqref{toeplitz} that for all $(x',\xi')$ in $S^*\partial X$
$$\sup_{\xi_n}|p_0(x',0,\xi',\xi_n) - f(x',\xi')|= \inf_{C\in K}\|p_0(x',0,\xi',D_n)-C -f(x',\xi')\mbox{Id}\|<\varepsilon.$$ 
Letting $\xi_n\to\infty$, the zero-homogeneity of $p_0$ implies that
the left hand side is $\ge |p_0(x',0,0,1)-f(x',\xi')|$. 
As $\varepsilon$ was arbitrary, 
$f$ is actually independent of the covariable: $f\in C(\partial X)$. This implies that $\lambda(F)$ belongs to $\mbox{im}\,b$.  Hence equality holds in (\ref{img}) and  
the image of $\gamma$ fits into the following exact sequence of  $C^*$-algebras
 \begin{equation}\label{eq:split_ex_seq}
   0\to C(S^*\partial X,{\mathcal W}_0)\to 
\mbox{im}\,\gamma  \to
   C(\partial X)\to 0.
 \end{equation}
This sequence splits via $b$. 
Now  Lemma \ref{kw} and the Künneth formula show that the $K$-theory of
$C_0(S^*\partial X,\mathcal{W}_0)$ vanishes, and 
Theorem~\ref{thmimg} follows from (\ref{eq:split_ex_seq}).

\section{K-Theory and index of Boutet de Monvel operators I}

We start this section recalling some results concerning the $K$-theory of
C$^*$-algebras, see \cite[Section 2]{MSS}. Let $\mathbb{A}$ be a $C^*$-algebra. The cone over $\mathbb{A}$ is the $C^*$-algebra 
$C\mathbb{A}:=\{\phi\colon [0,1]\to \mathbb{A}$; $\phi$ is continuous and $\phi(1)=0\}$. Since $C\mathbb{A}$ is a contractible $C^*$-algebra, its $K$-theory vanishes. The suspension of $\mathbb{A}$ is given by     
$S\mathbb{A}:= \{\phi\in C\mathbb{A}; \phi(0)=0\}$. If $f\colon \mathbb{B}\to \mathbb{A}$ is a $C^*$-algebra homomorphism, the mapping cone $Cf$ is defined to be 
$Cf:=\{ (b,\phi)\in \mathbb{B}\oplus C\mathbb{A};\;f(b)=\phi(0)\}$. The projection $q$ onto $B$ defines a short exact sequence 
\begin{equation}\label{eq:cone_seq}
0\,\longrightarrow\, S\mathbb{A}\,{\mathop{\longrightarrow}\limits^{i}}\, Cf{\mathop{\longrightarrow}\limits^{q}}\,\mathbb{B}\,\longrightarrow 0, 
\end{equation}
with $i$ denoting the inclusion $i:S\mathbb{A}\ni\phi\mapsto (0,\phi)\in Cf$. 

The assignment of the exact sequence (\ref{eq:cone_seq}) to each C$^*$-algebra
homomorphism $f\colon \mathbb{B}\to\mathbb{A}$ defines a functor between the corresponding categories (whose morphisms consist of commutative diagrams of homomorphisms or of exact sequences, respectively). This functor is exact. Another important observation 
is the following:

\begin{lem}\label{connmap}
The connecting maps in the standard cyclic 6-term exact sequence associated to \eqref{eq:cone_seq} are equal, modulo the canonical isomorphisms $K_i(S\mathbb{A})\cong K_{1-i}(\mathbb{A})$, to the group homomorphisms induced by $f$. 

If $f$ is additionally surjective, then the map 
$j\colon \ker f \to Cf$, given by $x\mapsto (x,0)$,
induces a $K$-theory isomorphism, which fits into the commutative diagram
$$
\def\mapup#1{\Big\uparrow\rlap{$\vcenter{\hbox{$\scriptstyle#1$}}$}}
\begin{array}{ccccccccc}
\longrightarrow&K_{i+1}(B)&\longrightarrow&K_i(SA)&{\mathop{\longrightarrow}\limits^{\hat\iota_*}}&K_i(Cf)&
{\mathop{\longrightarrow}\limits^{q_*}}&K_i(B)&\longrightarrow\\
                &\mapup{=}&     &\mapup{\delta_{i+1}}&          &\mapup{j_*}&             &\mapup{=}&\\
\longrightarrow&K_{i+1}(B)&{\mathop{\longrightarrow}\limits^{f_*}}&K_{i+1}(A)&\longrightarrow&K_i(\ker f)&\longrightarrow&K_i(B)&\longrightarrow
\end{array}
$$
where the upper row is the cyclic exact sequence induced by \eqref{eq:cone_seq}, and the lower one is that induced by
$
0\longrightarrow \ker f\longrightarrow B{\mathop{\longrightarrow}\limits^{f}} A\longrightarrow 0$.
\end{lem}

Applying the short exact sequence \eqref{eq:cone_seq} to the commutative diagram 
\eqref{commdiag}, one obtains the commutative grid:
 \begin{equation}\label{cde}
    \begin{CD}
      && 0 && 0 && 0\\
      && @AAA @AAA @AAA\\
      0 @>>> C_0(X^\circ) @>>> C(X) @>r>> C(\partial X) @>>> 0\\
      && @AAA @AAA @AAA\\
      0 @>>> Cm^\circ @>>> Cm @>>> Cb @>>> 0\\
      && @AAA @AAA @AAA\\
      0 @>>> S(\is/\ks) @>>> S(\A/\ks) @>{S\pi}>>
      S(\A/\is) @>>> 0\\
      && @AAA @AAA @AAA\\
      && 0 && 0 && 0
    \end{CD}\ .
  \end{equation}
Next consider the associated long exact sequences in $K$-theory.
By Theorem \ref{thmimg}, $b$ induces an isomorphism in $K$-theory. 
>From Lemma~\ref{connmap} and the cyclic exact sequence of 
$0\to S(\A/\is)\to Cb\to C(\partial X)\to 0$ 
we then conclude that $K_*(Cb)=0$. From this in turn we
deduce, using the cyclic exact sequence of $0\to Cm^\circ \to Cm\to
Cb\to 0$, that $Cm^\circ\to Cm$ induces an isomorphism in $K$-theory.

We know from Theorem \ref{kernelgamma} 
that $\is/\ks=\ker \gamma/\ks\cong C_0(S^*X^\circ)$.  Together with the canonical isomorphism
$K_*(S(\is/\ks))\cong K_{1-*}(\is/\ks)$ the left two vertical exact sequences 
induce the following commutative diagram in $K$-theory:
\begin{equation}
\label{2ces}
\def\mapdown#1{\downarrow\rlap{$\vcenter{\hbox{$\scriptstyle#1$}}$}}
\def\mapup#1{\uparrow\rlap{$\vcenter{\hbox{$\scriptstyle#1$}}$}}
\begin{array}{ccc}
K_0(C_0(X^\circ))          &\longrightarrow                         &K_0(C(X))\\
\mapdown{m^\circ_*}&                                        &\mapdown{m_*}\\
K_0(C_0(S^*X^\circ))     &\longrightarrow &K_0(\cA/\ks)\\
{\downarrow}          &                                        &{\downarrow}\\
K_1(Cm^\circ)              &{\mathop{\longrightarrow}\limits^{\cong}}&K_1(Cm)\\
\downarrow               &                                        &\downarrow\\
K_1(C_0(X^\circ))          &\longrightarrow                         &K_1(C(X))\\
\mapdown{m^\circ_*}&                                        &\mapdown{m^\circ_*}\\
K_1(C_0(S^*X^\circ))     &\mathop{\longrightarrow}\limits^{i_*}&K_1(\mathcal{A}/\mathcal{K})\\
 \mapdown{\alpha}         &                                        &\mapdown{\beta}\\
K_0(Cm^\circ)              &{\mathop{\longrightarrow}\limits^{\cong}_\phi}&K_0(Cm)\\
\downarrow               &                                        &\downarrow\\
K_0(C_0(X^\circ))          &\longrightarrow                         &K_0(C(X)).
\end{array}
\end{equation}

We shall now see how this can be used to derive a $K$-theoretic proof of 
Boutet de Monvel's index theorem. 
A crucial ingredient is the following well-known result:

\begin{lem}\label{section} 
A connected compact manifold with non-empty boundary always has a nowhere vanishing vector field. 
\end{lem} 

This implies that the co-sphere bundle of $X$ has a continuous section
$s$. Composition with $s$ then defines a left inverse for $m^\circ \colon C_0(X^\circ)\to C_0(S^*X^\circ)$. 
This yields a right inverse $s'$ for the map $\alpha$ in \eqref{2ces}, and $s''=i_*\circ s'\circ \phi^{-1}$ yields a right inverse for $\beta$.
Hence both long exact sequences in \eqref{2ces} split, and that on the right hand side 
yields the split short exact sequences 
\begin{eqnarray}\label{ses}
0\,\longrightarrow\, K_i(C(X))\,{\mathop{\longrightarrow}\limits^{m_*}}\, K_i(\cA/\ks){\mathop{\longrightarrow}\limits^{\beta}}\,K_{1-i}(Cm)\,\longrightarrow 0,\ \ i=0,1. 
\end{eqnarray}
It is worthwhile noting an immediate consequence of this split exactness:
\begin{cor}\label{sum}Each element of $K_i(\cA/\ks)$ can be written as the sum of an element 
in the range of $m_*$ and an element in the range of $s''$, hence in the range of 
$i_*$. 
\end{cor}

In order to determine 
$K_i(Cm)\cong K_i(Cm^\circ)$ we consider the commutative diagram 
\begin{equation}\label{tbs}
  \begin{CD}
    0 @>>>  C_0(T^*X^\circ) @>>> C_0(B^*X^\circ) @>{r}>>
    C_0(S^*X^\circ) @>>> 0\\
    && && @A{\pi^* r_0}A{\sim}A @AA{=}A\\
    && && C_0(B^*X^\circ) @>{\pi^* r_0}>> C_0(S^*X^\circ)\\
    && && @V{r_0}V{\sim}V @VV{=}V\\
    && && C_0(X^\circ) @>{\pi^*=m^\circ} 
                                   >> C_0(S^*X^\circ).
\end{CD}
\end{equation}
Here, $\pi^*$ denotes pull back from the base to the total space of
the bundle, while $r$ and $r_0$ denote restriction to the boundary of the ball
bundle and the zero section of the ball bundle, respectively; $\sim$ denotes homotopy
equivalence of $C^*$-algebras. 

We get induced short exact mapping cone sequences
\begin{equation}\label{eq:grid}
  \begin{CD}
    0 @>>> SC_0(S^*X^\circ) @>>>  Cr @>>> C_0(B^*X^\circ) @>>> 0\\
    && @AA{=}A @AA{(\pi^*r_0)_*}A @A{\sim}A{\pi^*r_0}A \\
  0 @>>> SC_0(S^*X^\circ) @>>>   C(\pi^* r_0) @>>> C_0(B^*X^\circ)
  @>>> 0\\
  && @VV{=}V  @V{\sim}V{(r_0)_*}V @V{\sim}V{r_0}V\\
 0 @>>> SC_0(S^*X^\circ) @>>> Cm^\circ  
                                    @>>> C_0(X^\circ) @>>> 0.
  \end{CD}
\end{equation}

Applying the 5-lemma to the corresponding cyclic exact $K$-theory sequences  we see
that the induced maps between the mapping cones are
$K$-theory isomorphisms. 

Finally, since $r$ is surjective and $\ker{r}=C_0(T^*X^\circ)$, Lemma~\ref{connmap} yields the commutative digram
\begin{equation}\label{cdj}
\def\mapup#1{\Big\uparrow\rlap{$\vcenter{\hbox{$\scriptstyle#1$}}$}}
\begin{array}{ccc}
K_0(SC_0(S^*X^\circ)) &\longrightarrow&K_0(Cr)\\
{\iso\Big\uparrow} &&\iso\mapup{j_*}\\
K_1(C_0(S^*X^\circ))&{\mathop{\longrightarrow}\limits^{\delta}}&K_0(C_0(T^*X^\circ))
\end{array},
\end{equation}
where the lower horizontal arrow is the index mapping for the first row in \eqref{tbs}, and the upper horizontal is induced by 
the first row in \eqref{eq:grid}.

This furnishes natural isomorphisms
$$K_i(Cm)\cong K_i(Cm^\circ)\cong K_i(C_0(T^*X^\circ)).$$

We next consider the commutative diagram
\begin{equation}\label{index}
\def\mapdown#1{\Big\downarrow\rlap{$\vcenter{\hbox{$\scriptstyle#1$}}$}}
\def\mapup#1{\Big\uparrow\rlap{$\vcenter{\hbox{$\scriptstyle#1$}}$}}
\begin{array}{ccccccl}
 &K_1(C(X))&{\mathop{\longrightarrow}\limits^{m_*}}           &K_1(\A/\ks)&\xrightarrow{\beta}&K_0(Cm)         &\\
      &\mapup{} &                                             & \mapup{i_*}           &   &\mapup{\cong}    &\\
 &K_1(C_0(
X^\circ ))&{\mathop{\longrightarrow}\limits^{ {m_{*}^\circ}}}&K_1(\is/\ks) &\xrightarrow{\alpha}&K_0(C m^\circ)   & \\
      &                &                              &a\mapup{\iso}          &   &\mapup{\iso}    &\\
       &                &                             &K_0(SC_0(S^* X^\circ ))&\longrightarrow&K_0(C(\pi^*r_0))&\\
        &                &                            &\mapdown{=}             &   &\mapdown{\iso}  &\\
         &                &                           &K_0(SC_0(S^* X^\circ ))&\longrightarrow&K_0(Cr)  &       \\
          &                &                          &c\mapup{\iso}          &   &\mapup{\iso}    &\\
           &                &                         &K_1(C_0(S^* X^\circ ))   &{\mathop{\longrightarrow}\limits^{\delta}}&K_0(C_0(T^* X^\circ ))&\\
            &                &                        &                      &   &\mapdown{\ind_t}&\\
             &                &                       &                      &   &\mathbb Z&
\end{array}
\end{equation}                                   
where the first two rows are portions of \eqref{2ces}. The second, third and fourth rows in \eqref{index}\ are portions 
of the cyclic sequences associated to \eqref{eq:grid}\ (notice that, if we use the isomorphism $\is/\ks\iso
C_0(S^* X^\circ )$\ as an identification, then the first column in \eqref{cde}\ is equal to the last row in \eqref{eq:grid}),
while the fourth and fifth rows are just \eqref{cdj}. Note that the
composed
isomorphism $c^{-1}a^{-1}\colon K_1(\is/\ks)\to 
K_1(C_0(S^*X^\circ))$ in the left row is exactly the map induced by the
interior symbol.

\begin{df} We define the map $p\colon K_1(\A/\is)\to K_0(C_0(T^* X^\circ ))$
as the composition of all the maps
(reverting arrows of isomorphisms when necessary) in the right  column in 
\eqref{index}, except $\ind_t$, with the map $\beta$ from
  $K_1(\A/\ks)$\ to $K_0(Cm)$\ in the first row.
\end{df}

We then infer from \eqref{ses}:
\begin{thm} $K_i(\A/\ks)$ fits into the short exact sequence
\begin{equation}
0\,\longrightarrow\, K_i(C(X))\,{\mathop{\longrightarrow}\limits^{m_*}}\, K_i(\cA/\ks){\mathop{\longrightarrow}\limits^{p}}\,K_{1-i}(C_0(T^*X^\circ))\,\longrightarrow 0,\ \ i=0,1. 
\end{equation}
The sequence splits, but not naturally.
\end{thm}

For $i=1$, we thus have a natural map
\begin{equation}\label{kth}
K_1(\A/\ks)\mathop{\longrightarrow}\limits^{p} K(C_0(T^*X^\circ))
\cong K_c(T^*X^\circ),
\end{equation}
where the last isomorphism is the identification of C$^*$-algebra $K$-theory with 
compactly supported $K$-theory of topological spaces. We can now state:

\begin{thm}\label{Bindex}
Let $\chi: K_c(T^*X^\circ)\to \mathbb Z$ be the topological index map defined by 
Atiyah and Singer. For an elliptic boundary value problem $A\in \A$ we then have 
\begin{eqnarray}\label{Bindexformula}
 \ind A = \chi\circ p([[A]]_1).
\end{eqnarray}
\end{thm}
Here $[[A]]_1$ is the $K_1$-class of the class $[A]$ of $A$ in $\A/\ks$, 
and we have used the 
identification of the $K$-theories mentioned above.

\begin{rem}\rm Further analysis shows that this map is precisely the map 
Boutet de Monvel constructed in \cite{B}  using deformations of boundary value 
problems and topological $K$-theory. See  \cite[Section 4]{MSS} for details.
\end{rem}

In order to prove Theorem \ref{Bindex}, we note that, by Corollary \ref{sum}, it is 
sufficient to treat the two cases where $[A]$ is in the range of $m_*$ or in the range 
of $i_*$. The elements in the range of $m_*$ are equivalence classes of 
invertible multiplication operators. 
Their Fredholm index therefore is zero.  
On the other hand,  the first row in (\ref{index}) is exact, 
thus the range of $m_*$ is mapped to zero. Hence both sides of \eqref{Bindexformula}
are zero. 

If $[A]$ is in the range of the map $i_*$ induced by the inclusion 
$i\colon \is/\ks\hookrightarrow \A/\ks$, then we may assume that $A$ is of the form \eqref{defis}, 
and the equality of both sides in \eqref{Bindexformula} essentially follows from the 
Atiyay-Singer index theorem by considering $\varphi P \varphi$ as a pseudodifferential
operator on $\widetilde X$. 
 
This proof breaks down in the 
case of elliptic families. Then it will no longer be true that the map $m^\circ$ has a left 
inverse.  In the next section we will outline an alternative way of computing 
the $K$-theory of $\A/\ks$. This approach will extend to the families case and lead to a proof of an index theorem for families, as explained in Section~\ref{families}. 

\section{K-Theory and index of Boutet de Monvel operators II}

Let $B$ denote the subalgebra of $C(S^*X)$ consisting of the functions which
do not depend on the co-variable over the boundary, that is, $f\in C(S^*X)$
belongs to $B$ if and only if the restriction of $f$ to $S^*X|_{\partial X}$
is of the form $g\circ\pi$, for some $g\in C(\partial X)$, where $\pi\colon S^*X\to X$ is the canonical projection. We will denote by $S^*X/\!\!\sim$ the quotient of $S^*X$ by the equivalence relation which identifies all $e$, $f\in S^*X|_{\partial X}$ such that $\pi(e)=\pi(f)$. The algebra $B$ is then canonically isomorphic to $C(S^*X/\asim)$.

Let $\cA^\dagger$ denote the C$^*$-subalgebra of $\lh$ generated by all operators of the form
\begin{eqnarray}\label{punhal}
A=
\begin{pmatrix}P_+ &K_1\\K_2&K_3
\end{pmatrix}: 
\begin{array}{ccc}
L^2(X)& & L^2(X)\\
\oplus&\to&\oplus\\
L^2(\partial X)&&L^2(\partial X)
\end{array} ,
\end{eqnarray} 
where $P$ is a pseudodifferential operator satisfying the transmission condition with principal symbol belonging to $B$, and $K_1$, $K_2$ and $K_3$ are compact operators. Comparing with (\ref{defis}) it is then clear that 
$\is=\ker\gamma$ is contained in $\cA^\dagger$. 

\begin{pro}\label{hoc} The restriction of $\sigma$ to $\cA^\dagger$ has kernel equal to $\ks$ and image equal to $B$. In other words, the principal symbol 
induces an isomorphism $\cA^\dagger/\ks\cong B$. 
\end{pro}
{\em Proof}: If the upper left corner of an $A\in\cA^\dagger$ is $P_+$, where $P$ is a pseudodifferential operator satisfying the transmission condition with principal symbol belonging to $B$, then 
$\gamma(A)=\sigma(A)|_{S^*X|_{\partial X}}\otimes\mbox{Id}$, where $\mbox{Id}$ denotes the identity on 
$L^2(\R_{\geq 0})\oplus\C$. Since the set of all such $A$ generate $\cA^\dagger$ and $\gamma$ and $\sigma$ are homomorphisms, we have $\gamma(A)=\sigma(A)|_{S^*X|_{\partial X}}\otimes\mbox{Id}$ for all $A\in\cA^\dagger$. In particular, the kernel of the restriction of $\sigma$ to $\cA^\dagger$ is contained in $\ker\sigma\cap\ker\gamma=\ks$. It is equal to $\ks$ because it contains all integral operators with smooth kernel; they are Boutet de Monvel operators of order $-\infty$ and class zero. 

If $q\in B\cap C^\infty(S^*X)$, then $q$ is the principal symbol of a pseudodifferential operator satisfying the transmission condition (see \cite[Theorem~1 of Section~2.3.3.1]{RS}, for example). The algebra 
$B\cap C^\infty(S^*X)$ separates points in $S^*X/\asim$. Hence $\sigma(\cA^\dagger)$ is a dense subalgebra of $B$, which is a closed subalgebra of $C(S^*X)$. This finishes the proof, since the image of a C$^*$-algebra homomorphism is always closed. \cqd

The following proposition can be proven by a diagram chase (see \cite[Exercise 38, Section~2.2]{H}). 
\begin{pro}\label{ditsche} Let there be given a commutative diagram of 
abelian groups with exact rows, 
\[
\def\mapdown#1{\downarrow\rlap{$\vcenter{\hbox{$\scriptstyle#1$}}$}}
\def\mapup#1{\uparrow\rlap{$\vcenter{\hbox{$\scriptstyle#1$}}$}}
\begin{array}{ccccccccccc}
\cdots&\rightarrow
&A_i^\prime&{\mathop{\longrightarrow}\limits^{f_i^\prime}}
&B_i^\prime&{\mathop{\longrightarrow}\limits^{g_i^\prime}}
&C_i^\prime&{\mathop{\longrightarrow}\limits^{h_i^\prime}}
&A_{i+1}^\prime
&\rightarrow&\cdots
\\
&
&\mapup{a_i}&&\mapup{b_i}&&\mapup{c_i}&&\mapup{a_{i+1}}&
&
\\
\cdots&\rightarrow
&A_i&{\mathop{\longrightarrow}\limits^{f_i}}
&B_i&{\mathop{\longrightarrow}\limits^{g_i}}
&C_i&{\mathop{\longrightarrow}\limits^{h_i}}
&A_{i+1}
&\rightarrow&\cdots
\end{array},
\]
where each $c_i$ is an isomorphism. Then the sequence
\[
\cdots\longrightarrow A_i{\mathop{\longrightarrow}\limits^{(a_i,-f_i)}}
A_i^\prime\oplus B_i{\mathop{\longrightarrow}\limits^{\langle f_i^\prime,b_i\rangle}}
B_i^\prime{\mathop{\longrightarrow}\limits^{h_ic_i^{-1}g_i^\prime}}
A_{i+1}\longrightarrow\cdots
\]
is exact, where $\langle f_i^\prime,b_i\rangle$ is the map defined by
$\langle f_i^\prime,b_i\rangle(\alpha,\beta)=f_i^\prime(\alpha)+b_i(\beta)$.
\end{pro}

\begin{thm}\label{kthth}
Let $\iota:\cA^\dagger/\ks\to\cA/\ks$ denote the canonical inclusion. Then 
\[
\iota_*:K_*(\cA^\dagger/\ks)\to K_*(\cA/\ks) 
\]
is an isomorphism.
\end{thm}
In view of Proposition \ref{hoc} this furnishes a description of the 
$K$-theory of $\A/\ks$ in terms of that of a topological space.

{\em Proof}: 
Applying Proposition \ref{ditsche} to the diagram (\ref{2ces}), we get the exact sequence
\begin{equation}
\label{cyclic1}
\begin{array}{ccccc}
K_0(C_0(X^\circ))&\rightarrow &K_0(C(X))\oplus K_0(C_0(S^*X^\circ))&\rightarrow& 
K_0(\cA/\ks)
\\\uparrow& & & &\downarrow
\\
K_1(\cA/\ks)&\leftarrow& K_1(C(X))\oplus K_1(C_0(S^*X^\circ))&\leftarrow&  
K_1(C_0(X^\circ))
\end{array}.
\end{equation}

We next consider the following diagram of commutative C$^*$-algebras
\begin{equation}\label{cd3}
\def\mapdown#1{\downarrow\rlap{$\vcenter{\hbox{$\scriptstyle#1$}}$}}
\def\mapup#1{\uparrow\rlap{$\vcenter{\hbox{$\scriptstyle#1$}}$}}
\begin{array}{ccc}
C_0(X^\circ)&{\mathop{\longrightarrow}\limits^{m^\circ}}&C_0(S^*X^\circ)\\
\downarrow&&\mapdown{p_2}\\
C(X)&{\mathop{\longrightarrow}\limits^{p_1}}&B
\end{array}.
\end{equation}
As $C_0(X^\circ)$ is canonically isomorphic to
\[
\{(f,g)\in C(X)\oplus C_0(S^*X^\circ);\ p_1(f)=p_2(g)\},
\]
the Mayer-Vietoris exact sequence  associated to (\ref{cd3}) is the 
exact sequence
\begin{equation}
\label{cyclic2}
\begin{array}{ccccc}
K_0(C_0(X^\circ))&\rightarrow &K_0(C(X))\oplus K_0(C_0(S^*X^\circ))&\rightarrow& 
K_0(B)
\\\uparrow& & & &\downarrow
\\
K_1(B)&\leftarrow& K_1(C(X))\oplus K_1(C_0(S^*X^\circ))&\leftarrow&  
K_1(C_0(X^\circ))
\end{array}.
\end{equation}
The maps $\iota_*\colon K_i(B)\cong K_i(\cA^\dagger/\ks)\rightarrow K_i(\cA/\ks)$, $i=0,1$, and the identity on 
the other $K$-theory groups furnish morphisms from the 
cyclic sequence (\ref{cyclic2}) to  the cyclic sequence (\ref{cyclic1}).
The five lemma then shows that $\iota_*$ is an isomorphism. 
\cqd

Let $\beta\colon K_1(C(S^*X/\asim))\to K_0(C_0(T^*X^\circ))$ denote 
the connecting map in the standard cyclic exact sequence associated to 
\begin{equation}
\label{tbstil}
0\longrightarrow C_0(T^*X^\circ)\longrightarrow C(B^*X/\asim)
\longrightarrow C(S^*X/\asim)\longrightarrow 0,
\end{equation}
where $B^*X$ denotes the bundle of closed co-balls over $X$ (which can be regarded as a compactification 
of $T^*X$ whose points at infinity form the co-sphere bundle $S^*X$). 

Let $e\colon C_0(T^*X^\circ)\to C_0(T^*\widetilde X)$ be the map
of extension by zero and denote with
$j\colon K_0(C_0(T^*\widetilde X))\to K(T^*\widetilde X)$ the canonical isomorphism 
between C$^*$-algebra $K$-theory and topological $K$-theory groups. 

\begin{thm}\label{indthm} If $A\in\cA$ is a Fredholm operator, then
\begin{equation}
\label{indexformula}
\mbox{{\em ind}}\,A=\chi\circ j\circ e_*\circ\beta\circ\iota_*^{-1}([[A]]_1)
\end{equation}
where $\mbox{{\em ind}}$ denotes the Fredholm index, $[[A]]_1$ denotes the $K_1$ of the class 
$[A]$ of $A$ in the quotient $\cA/\ks$ and 
$\chi:K(T^*\widetilde X)\to \Z$ is Atiyah and Singer's topological index for the closed manifold 
$\widetilde X$. 
\end{thm}

In Theorem \ref{Bindex} above,  we stated that $\ind A = \chi\circ p( [[A]]_1)$,
using the identification of C$^*$-algebra $K$-theory and topological $K$-theory. 
Actually, we could have been more precise, because the 
Atiyah-Singer topological index map 
is only defined for the cotangent bundle of a {\em closed} manifold, and 
the identification involves the maps $e_*$ and $j$. Thus that formula 
should actually read 
\begin{eqnarray}\label{Bindexformula2}
\ind A = \chi\circ j\circ e_*\circ p ([[A]]_1).
\end{eqnarray}
This is what was in fact shown in \cite[Theorem 2]{MSS}.

We can infer Theorem \ref{indthm} from Equation \eqref{Bindexformula2} by 
showing that 
\begin{equation}\label{equality}
p=\beta\circ\iota_*^{-1}.
\end{equation}

For that, let us consider the following commutative diagram of exact sequences
\begin{equation}
\label{commdiag2}
\def\mapup#1{\Big\uparrow\rlap{$\vcenter{\hbox{$\scriptstyle#1$}}$}}
\begin{array}{ccccccc}
0\longrightarrow&C_0(T^*X^\circ)&\longrightarrow&C(B^*X/\asim)&{\mathop{\longrightarrow}\limits^{r}}
&C(S^*X/\asim)&\longrightarrow 0
\\                &\mapup{=}&&\mapup{i^b}&&\mapup{i^s}&
\\0\longrightarrow&C_0(T^*X^\circ)&\longrightarrow&C(B^*X^\circ)&{\mathop{\longrightarrow}\limits^{r}}
&C(S^*X^\circ)&\longrightarrow 0
\end{array}.   
\end{equation}
Recall that $\beta$  and $\delta$ are the index maps associated to the upper and 
lower sequence, respectively. 
The naturality of the index map implies that 
\[
\delta=\beta\circ i^s_*.
\] 
Up to the isomorphisms of Theorem~\ref{kernelgamma} and Proposition~\ref{hoc}, the map $i^s$ in (\ref{commdiag2}) is equal to the canonical inclusion $\is/\ks\to\cA^\dagger/\ks$. 
With $i\colon\is/\ks\hookrightarrow \cA/\ks$ we have 
\[
p\circ i_*=\delta
\]
up to the $K$-theory isomorphism induced by the C$^*$-algebra isomorphism of Theorem~\ref{kernelgamma}: For that, see the diagram \eqref{index} and the remark following it. 
Since $i=\iota\circ i^s$, this shows that 
\[
p\circ\iota_*\circ i^s_*=\beta\circ i_*^s;
\]
that is, $p\circ\iota_*=\beta$ holds on the image of $i^s_*$; or, equivalently, (\ref{equality}) holds on the image of $i_*$. 

In view of Corollary \ref{sum}, it remains to show that $p\circ \iota_*=\beta$ also 
holds on the image of $m_*$.
Now, the exactness of (\ref{kth}) implies that $p\circ m_*=0$. Hence, all that is left to 
prove Theorem~\ref{indthm} is to show that $\beta\circ m_*=0$. If an element $x\in K_1(C(S^*X/\asim))$ is represented by an invertible $f\in M_k(C(X))$ (notice that we are using the isomorphism of 
Proposition~\ref{hoc} as an identification), then $x$ belongs to the image of $r_*:K_1(C(B^*X/\asim))\to K_1(C(S^*X/\asim))$ (since $f$ can also be regarded as an invertible of $M_k(C(B^*X/\asim))$). The exactness of the cyclic exact sequence associated to (\ref{tbstil}) implies that $\beta(x)=0$.

Theorem~\ref{indthm} is also a particular case (when the space of parameters $Y$ reduces to one point) of 
Theorem~\ref{famindthm}, which can be proven independently of 
Theorem \ref{Bindex}.

\section{K-theory and index for families}\label{families}

The index of a continuous function taking values in Fredholm operators acting on a fixed Hilbert space was defined by J\"anich \cite{J} and Atiyah \cite{A}. That definition was adapted by Atiyah and Singer \cite{AS4} to continuous families of elliptic operators acting on the fibers of a fiber bundle whose fibers are closed manifolds. A slight variation of their definitions, for sections of Fredholm operators in a bundle of 
C$^*$-algebras, is used in \cite{MSSf} to state and prove Theorem~\ref{famindthm} below, which is based on and generalizes Atiyah and Singer's index theorem for families.

Let $X$, as before, be a compact manifold with boundary, 
and take $\widetilde X=2X$, the double of $X$. Let $Y$ be a compact Hausdorff space and let 
$\pi\colon Z\to Y$ be a fiber bundle with fiber $X$ and structure group $\mbox{Diff}(X)$ (equipped with its usual topology). Each $Z_y=\pi^{-1}(y)$ is a compact manifold with boundary, noncanonically diffeomorphic to $X$. Let $\delta\colon U\to\partial X\times [0,1)$ be a diffeomorphism defined on an open neighborhood of $\partial X$. The structure group of the bundle $\pi$ can be reduced \cite[Appendix A]{MSSf} to the subgroup $G$ of $\diff(X)$ consisting of all those $\phi$ such that $\delta\circ\phi\circ\delta^{-1}\colon\partial X\times[0,1/2)\to\partial X\times[0,1)$ is of the form $(x^\prime,x_n)\mapsto (\varphi(x^\prime),x_n)$ for some diffeomorphism 
$\varphi\colon\partial X\to \partial X$. The elements of $G$ are such that their reflections to the double $2X$ of $X$ are also diffeomorphisms and this allows us to consider the double $2Z$ of $Z$, a fiber bundle over $Y$ with fiber $2X$ and structure group $\mbox{Diff}(2X)$. Another consequence of this technicality is that we will then be able to define the boundary principal symbol of a family of Boutet de Monvel (as remarked after 
(\ref{bdysymbol}), the boundary principal symbol is invariantly defined only after we choose a normal coordinate $x_n$). 

We next fix a continuous family of riemannian metrics on $Z_y$ and use them to define the Hilbert spaces 
$H_y=L^2(Z_y)$. 
The union ${\displaystyle \hc=\bigcup_{y\in Y}H_y}$ can be canonically given the 
structure of a fiber bundle with fiber $H=L^2(X)$ and structure group $G$; here 
$G$ acts on $L^2(X)$ via the representation $\phi\mapsto T_\phi$, 
$T_\phi(f)=f\circ\phi^{-1}$, which is continuous with respect to the strong operator 
topology. 

To simplify the exposition, we will denote, for the rest of this section, by $\A$ only the upper left corner of what was denoted by $\A$ in Definition~\ref{defas}; i.e., $\A$ is the norm closure in $\mathcal{L}(L^2(X))$ of the algebra of all operators $P_++G$, where $P$ is a classical pseudodifferential operator of order zero satisfying the transmission condition and $G$ is a polyhomogeneous singular Green operator of order and class zero. Analogously, for each $y\in Y$, we define $\A_y$ as the norm closure in the bounded operators on $L^2(Z_y)$ of (the upper left corner of) the algebra of all Boutet de Monvel operators on $Z_y$.

It is well known that the Boutet de Monvel algebra is invariant under diffeomorphisms. Arguing similarly as in \cite[Proposition 1.3]{AS4}, one can show, furthermore, that the map
\begin{eqnarray}\label{jointcont}
G\times \A\ni (\phi,A) \mapsto T_\phi A T^{-1}_\phi \in\A
\end{eqnarray}
is jointly continuous. This implies that the union ${\displaystyle \aleph=\bigcup_{y\in Y}\A_y}$ can be canonically given the structure of a fiber bundle with fiber $\A$.

\begin{df}
The continuous sections of the bundle $\aleph$ form a C$^*$-algebra, which we denote by $\mathfrak{A}$.
\end{df}

\begin{rem}\rm Our approach differs slightly from that of Atiyah and Singer, who work 
with sections of Fréchet algebras instead of C$^*$-algebras.

Note that the continuity property \eqref{jointcont} is crucial and limits the 
choice of the algebras in the fibers. It is not possible, for example, to give 
$\bigcup_y \mathcal{L}(L^2(Z_y))$ in a canonical way the
structure a bundle of $C^*$-algebras (with structure group the unitary group  
with norm topology). 
\end{rem}

We recall the construction of the analytic index of families for the present situation.

Similarly as in \cite[Proposition (2.2)]{AS4} and \cite[Proposition A5]{A}, one can prove:
\begin{pro}\label{perturb} Let $\hc$ and $\ac$ be as above and let 
$(A_y)_{y\in Y}\in\ac$ be such that, for each $y$, 
$A_y$ is a Fredholm operator on $H_y$. Then there are continuous sections
$s_1,\cdots,s_q$ of $\mathfrak{H}$ such that the maps
\[
\begin{array}{rccl}
\tilde A_y:&H_y\oplus\C^q&\longrightarrow&H_y\oplus\C^q\\
              &(v,\lambda)&\longmapsto&(A_yv+\sum_{j=1}^{q}\lambda_js_j(y),0)
\end{array}
\]
have image equal to $H_y\oplus 0$ for all $y\in Y$ and hence
$(\ker\tilde A_y)_{y\in Y}$ is a finite-dimensional vector bundle 
over $Y$. 
\end{pro}

\begin{df}\label{ani} 
Given $A=(A_y)_{y\in Y}\in\ac$ as in Proposition~\ref{perturb}, we
denote by $\ker\tilde A$ the bundle $(\ker\tilde A_y)_{y\in Y}$ and
define
\[
\ind_a A=[\ker\tilde A]-[Y\times\C^q]\in K(Y).
\]
This is independent of the choices of $q$ and of $s_1,\cdots,s_q$ and 
we call it {\em the analytical index of} $A$.
\end{df}

If $k$ is an integer, the algebra $M_k(\ac)$ is naturally isomorphic to the 
algebra of continuous sections of the bundle of algebras 
$
{\displaystyle \aleph_k=\bigcup_{y\in Y}M_k(\mathcal{A}_y)}, 
$ 
each $M_k(\mathcal{A}_y)$ a C$^*$-subalgebra of the bounded operators
on $H_y^k$. We can then define $\ind_a(A)$ if $A=(A_y)_{y\in Y}\in M_k(\ac)$ 
is a section such that each $A_y$ is a Fredholm operator on $H_y^k$. 
The projection of such an $A=(A_y)_{y\in Y}\in M_k(\ac)$ in $M_k(\ac/\kc)$ 
is invertible and hence defines an element of $K_1(\ac/\kc)$. Since $\ind_a(A)$
is invariant under homotopies and pertubations by 
compact operator valued sections, we get a homomorphism 
\begin{equation}
\label{topind}
\ind_a\colon K_1(\ac/\kc)\longrightarrow K(Y).
\end{equation}

Let us denote by $S^*Z$ the disjoint union of all $S^*Z_y$. This can 
canonically be viewed as the total space of a fiber bundle over $Y$ with 
structure group $G$. 
One analogously defines $S^*\partial Z=\cup_{y}S^*\partial Z_y$ and 
$S^*Z^\circ=\cup_{y} S^*Z_y^\circ$. The families of homomorphisms 
\[
\sigma_y\colon \A_y\rightarrow C(S^*Z_y)
\ \ \mbox{and}\ \ 
\gamma_y\colon \A_y\rightarrow C(S^*\partial Z_y,\mathcal{L}(L^2(\R_{\ge0}))),\ \ y\in Y,
\]
piece together, yielding C$^*$-algebra homomorphisms 
\[
\sigma\colon \ac\longrightarrow C(S^*Z)
\ \ \mbox{and}\ \ 
\gamma\colon \ac\longrightarrow C(S^*\partial Z,\mathcal{L}(L^2(\R_{\ge0}))).
\]
For this, note in particular that $\gamma_y$  is well defined because the structure group of the bundle $\pi\colon Z\to Y$
leaves the normal coordinate of $X$ invariant, see \cite[Theorem 2.4.11]{G}.

Similarly as in equation \eqref{punhal} and Proposition \ref{hoc} 
we define $\ac^\dagger$ as the C$^*$-algebra generated by 
the families of pseudodifferential operators with principal symbol independent of the covariable over boundary points and show that the principal symbol 
$\sigma$ induces an isomorphism 
\[
\bar\sigma: \ac^\dagger/\kc\,\longrightarrow\,C(S^*Z/\asim),
\]
where $\kc$ denotes the continuous compact operator valued sections of $\aleph$ and 
$S^*Z/\asim$ denotes the union of all $S^*Z_y/\asim$, $y\in Y$.
The same arguments that prove Theorem~\ref{kthth} also prove that the canonical inclusion $\iota:\ac^\dagger/\kc\to\ac/\kc$ induces
a $K$-theory isomorphism 
\[
\iota_*:K_*(\ac^\dagger/\kc)\,{\mathop{\longrightarrow}\limits^{ }}\, K_*(\ac/\kc).
\]

Analogously as we did around (\ref{tbstil}), we also denote by 
$\beta\colon K_1(C(S^*Z/\asim))\to K_0(C_0(T^*Z^\circ))$ 
the index map in the standard cyclic exact sequence associated to 
\[
0\longrightarrow C_0(T^*Z^\circ)\longrightarrow C(B^*Z/\asim)
\longrightarrow C(S^*Z/\asim)\longrightarrow 0.
\]

\begin{df}\label{indt} 
The topological index {\em $\ind_t$} is the following composition of maps 
\begin{equation*}
\label{anaind}
\def\mapdown#1{\downarrow\rlap{$\vcenter{\hbox{$\scriptstyle#1$}}$}}
\def\mapup#1{\uparrow\rlap{$\vcenter{\hbox{$\scriptstyle#1$}}$}}
\begin{array}{rc}
\ind_t\colon K_1(\ac/\kc)
{\mathop{\longrightarrow}\limits^{\bar\sigma_*\circ\iota_{*}^{-1}}}
K_1(C(S^*Z/\asim )){\mathop{\longrightarrow}\limits^{\beta}} 
K_0(C_0(T^*Z^\circ)){\mathop{\longrightarrow}\limits^{e_*}}
&K_0(C_0(T^*2Z))\\ &
\mapdown{\mbox{{\sc as}}-\ind_t}\\ &K(Y),
\end{array}
\end{equation*}
where $e\colon C_0(T^*Z^\circ)\to C_0(T^*2Z)$ is the map
which extends by zero and {\em $\mbox{{\sc as}}-\ind_t$}
denotes the composition of Atiyah and Singer's \cite{AS4}\ topological families-index for 
the bundle of closed manifolds $2Z$  with the canonical 
isomorphism $K(T^*2Z)\simeq K_0(C_0(T^*2Z))$.
\end{df}

We are ready to state the main result of \cite{MSSf}:

\begin{thm}\label{famindthm} The two homomorphisms{\em 
\[ 
\ind_a:K_1(\ac/\kc)\to K(Y)\ \ \ \ \mbox{{\em and}}\ \  \ \ 
\ind_t:K_1(\ac/\kc)\to K(Y)
\]} 
are equal. \end{thm}

An arbitrary element of $K_1(\ac/\kc)$ is of the form $[[A]]_1$, where the inner brackets denote the class modulo compacts of a Fredholm operator valued element $A$ of $M_k(\ac)$. Our strategy to prove 
Theorem~\ref{famindthm} in \cite{MSSf} is to derive the equality of the
indices from the classical Atiyah-Singer index theorem for families
\cite[Theorem (3.1)]{AS4}.
To this end we defined a continuous family (in the sense of \cite{AS4}) 
of pseudodifferential operators $\hat A$ acting on a suitably constructed vector bundle over $2Z$ such that the topological indices of $A$ and of $\hat A$ are equal, and similarly the analytical indices of $A$ and $\hat A$ are also equal.

\vskip0.5cm

\tiny{
\noindent
Severino T. Melo. Instituto de  Matem\'atica e Estat\'{\i}stica, Universidade de S\~ao Paulo, Rua do Mat\~ao 1010, 05508-090 S\~ao Paulo, Brazil. E-mail: toscano@ime.usp.br.
\vskip0.2cm

\noindent
Thomas Schick. Mathematisches Institut, Universit{\"a}t G{\"o}ttingen, Bunsenstr.\ 3-5, 37073 G\"ottingen, Germany. E-mail: schick@uni-math.gwdg.de.
\vskip0.2cm

\noindent
Elmar Schrohe. Institut f\"ur Mathematik, Universit\"at Hannover, Welfengarten 1, 30167 Hannover, Germany. 
E-mail: schrohe@math.uni-hannover.de.}


\begin{thebibliography}{1}

\bibitem{A}{\sc M. F. Atiyah}. $K$-Theory, Lecture notes by D. W. Anderson. 
W. A. Benjamin, Inc., New York-Amsterdam, 1967. 

\bibitem{AS1} {\sc M. F. Atiyah \&\ I. M. Singer}. {\em The index of elliptic operators I}.
Ann. of Math. (2) {\bf 87} (1968), 484--530.

\bibitem{AS4}{\sc M. F. Atiyah \&\ I. M. Singer}. 
{\em  The index of elliptic operators IV}. 
Ann. of Math. (2)  {\bf 93} (1971), 119--138. 

\bibitem{Bl}{\sc B. Blackadar}. 
$K$-Theory for Operator Algebras. 
Cambridge University Press, Cambridge, 1998. 

\bibitem{B}{\sc L. Boutet de Monvel}. 
{\em Boundary problems for pseudo-differential operators}. 
Acta Math. {\bf 126} (1971), no. 1-2, 11--51.

\bibitem{Dou}{\sc R. G. Douglas}. Banach Algebra Techniques in Operator Theory, Grad. Texts Math. {\bf 179}, Springer, New York, 1998.

\bibitem{Ex} {\sc R. Exel}. {\em A Fredholm operator approach to Morita equivalence}. K-Theory {\bf 7}-3 (1993), 285-308.

\bibitem{F}{\sc B.~V.~Fedosov}. {\em Index theorems}.
In Partial differential equations VIII, 155--251, 
Encyclopaedia Math. Sci. {\bf  65}, Springer, Berlin, 1996.

\bibitem{Goh} {\sc I. C. Gohberg}. {\em On the theory of multidimensional singular integral operators}. Dokl. Akad. Nauk SSSR {\bf 133} (1960), 1279-1282.

\bibitem{G} {\sc G.~Grubb}. Functional Calculus of Pseudodifferential Boundary
Problems, Second Edition, Birkh{\"a}user, Boston, 1996.

\bibitem{GG} {\sc G. Grubb \&\ G. Geymonat}. {\em The essential spectrum of elliptic systems of mixed order}. Math. Ann. {\bf 227} (1977), 247-276.

\bibitem{H} {\sc A. Hatcher}. Algebraic Topology. Cambridge University Press, Cambridge, 2002. 

\bibitem{Hp} {\sc L. H\"ormander}. {\em Pseudo-differential \ops\ and hypoelliptic equations}. In Singular Integrals (Chicago, 1966), 
Proc. Symp. Pure Math. {\bf 10}, 138-183, 1967.

\bibitem{J}{\sc K. J\"anich}. 
{\em Vektorraumb\"undel und der Raum der Fredholm-Operatoren}. Math. Ann. {\bf 161} (1965), 129--142. 

\bibitem{KN} {\sc J. J. Kohn \&\ L. Nirenberg}. {\em An algebra of pseudo-differential operators}. Comm. Pure Appl. Math. {\bf 18} (1965), 269--305.

\bibitem{MNS}{\sc S. T. Melo, R. Nest, E. Schrohe}. 
{\em C$^*$-structure and $K$-theory of Boutet de Monvel's algebra}. 
J. reine angew. Math. {\bf 561} (2003), 145--175.

\bibitem{MSS}{\sc S. T. Melo, T. Schick, E. Schrohe}. 
{\em A $K$-theoretic proof of Boutet de Monvel's index theorem 
for boundary value problems}. 
J. reine angew. Math. {\bf 599} (2006), 217--233.

\bibitem{MSSf}{\sc S. T. Melo, T. Schick, E. Schrohe}. 
{\em Families index for Boutet de Monvel operators}. {\tt arXiv:1203.0482}.

\bibitem{RS} {\sc S.~Rempel and B.-W.~Schulze}. Index theory of elliptic boundary problems, Akademie-Verlag, Berlin, 1982.


\bibitem{S2} {\sc E.~Schrohe}. {\em Fr{\'e}chet algebra techniques for boundary value 
problems on noncompact manifolds: Fredholm criteria and functional 
calculus via spectral invariance}. Math. Nachr.  {\bf 199} (1999), 145--185.

\bibitem{S3} {\sc E.~Schrohe}.  {\em A short introduction to Boutet de Monvel's calculus}. 
In Approaches to Singular Analysis (Berlin, 1999),  Oper. Theory Adv. Appl. {\bf 125}, 85-116, Birkh{\"a}user, Basel, 2001.

\bibitem{See} {\sc R. T. Seeley}. {\em Integro-differential operators on vector bundles}. Trans. Amer. Math. Soc. {\bf 117} (1965), 167-204.

\bibitem{WO} {\sc N. E. Wegge-Olsen}. K-theory and C$^*$-algebras, Oxford University Press, New York 1993.

\end{thebibliography}
\end{document}